\documentclass[reqno]{amsart}
\author{Julia Brandes}
\thanks{The original version of this was submitted as a diploma thesis to the University of Stuttgart in 2009.}
\title{Twins of $s$-free numbers}

\address{School of Mathematics, University of Bristol, University Walk, Clifton, Bristol BS8 1TW, United Kingdom}
\email{mazjb@bristol.ac.uk}
\usepackage[T1]{fontenc}
\usepackage[latin1]{inputenc}
\usepackage{latexsym}
\usepackage[arrow, matrix, curve]{xy}
\usepackage[dvips]{graphicx}
\usepackage{epsfig}
\usepackage{bbm}
\usepackage{amssymb}
\usepackage{amsmath}
\usepackage{verbatim}
\usepackage{amsthm}
\usepackage{enumerate}

\def\res#1#2#3{ {\left( {\frac{#2}{#3}} \right) \negthickspace }_{#1} }
\def\dsum#1#2{ \sum_{\substack{ {#1} \\ {#2} }} }
\def\sP{{\mathcal{P}}}

\def\ksum{{\sum_{K<k \leq 2K}}}
\def\vsum{{\sum_{L \leq v \leq M}}}
\def\e#1#2{{}e\left( \frac{#1}{#2}\right) }

\newtheorem{satz}{Theorem}
\newtheorem{lem}{Lemma}

\theoremstyle{definition}

\newtheorem*{bems}{Remarks}
\newtheorem*{bem}{Remark}

\newenvironment{bew}{\begin{proof}[Proof]}{\end{proof}}

\begin{document}
\maketitle

\begin{abstract}
We generalise the square sieve developed by Heath-Brown to higher powers and apply it in order to improve on the error term for the problem of counting consecutive power-free numbers.
\end{abstract}

\section{Introduction}
Let $x>0$. How many numbers $n \leq x$ not divisible by any $s$-th power can there be? In the case $s=2$ the squared Möbius function acts quite naturally as an indicator and it does not require much effort to prove that
	$$\sum_{n \leq x} \mu^2(n)  =  \frac{6x}{\pi^2}  + O(\sqrt{x}).$$
This result can be generalized to larger $s$ and then yields
	$$ \# \left\lbrace n \leq x : \; n \hbox{ is } s \hbox{-free}  \right\rbrace = \frac{x}{\zeta(s)} + O(x^{1/s}),$$
where $\zeta(s)$ denotes the Riemann zeta function. But since the error term depends on the behaviour of $\zeta(s)$, it is impossible to achieve substantial improvements there without assuming the Riemann hypothesis to be true.
A synopsis of these and other results concerning $s$-free numbers is given in a survey paper by Pappalardi \cite{pappa}.

If, instead of counting $s$-free numbers one by one, one considers pairs with a given distance $a$, it turns out to be possible to improve on the result without having to assume the validity of the Riemann hypothesis. Let $E_s(n)$ be the indicator on the $s$-free numbers, then this means finding an estimate for
	$$\sum_{n\leq x} E_s(n) E_s(n+a).$$
In the most natural case $a=1$, an elementary result states that
\begin{align}
 	\sum_{n \leq x} E_s(n) E_s(n+1) = C_sx + O \big(x^{\frac{2}{s+1} + \epsilon} \big), \label{carlitz}
\end{align}
where the constant $C_s$ is given by the Euler product
\begin{align}
 	C_s = \prod_{p} \left( 1-\frac{2}{p^s} \right), \label{cs}
\end{align}
see Carlitz \cite{carl}.

In a paper from 1984, Roger Heath-Brown \cite{hb} succeeded to improve the error term in the case $s=2$ to  $O\left( x^{\frac{7}{11}}(\log x)^7 \right)$, using a square sieve. By similar methods we will show the following theorem.
\begin{satz}\label{satz}
 	For any $s \geq 2$ let $E_s(n)$ be the indicator function on $s$-free numbers, and $C_s$ the term given in \eqref{cs}. Then for any $\epsilon > 0$ one has
		$$\sum_{n \leq x} E_s(n) E_s(n+1) = C_sx + O \big( x^{{\frac{14}{7s+8}} + \epsilon} \big).$$
\end{satz}

Note that for $s=2$ this reproduces Heath-Brown's result. The exponent is $O(1/s^2)$ better than the one obtained by Carlitz, with $\epsilon$ accumulating the arising powers of $\log x$. These could also be expressed explicitly by calculating more precisely; in our context, however, this does not seem necessary, as the main saving occurs in the powers of $x$.\\

I would like to thank Prof. Dr. Jörg Brüdern for suggesting this problem to me.

\section{A sieve for $s$-th powers}

In the quadratic case Heath-Brown \cite{hb} applies a sieve based on Jacobi symbols. In order to generalize his method to higher powers, it is necessary to introduce an appropriate function that operates in a similar way on the $s$-th powers modulo $p$.

For $p$ with $(s, p-1) \neq 1$ let $\res snp$ be a non-principal character modulo $p$ such that $ \res s{n^s}p = \chi_0(n) $. Such a character exists, since the number of characters $\chi \pmod p$ with $\chi^s = \chi_0$ is given by $(p-1, s) \geq 2$.

Let $w$ be a nonnegative weight function on the integers such that $\sum_n w(n) < \infty$, and let $\mathcal{A}$ denote the sequence $\left(w(n) \right)$. We define
	$$S(\mathcal{A}) = \sum_{n=1}^{\infty} w(n^s).$$
With this notation, we have
\begin{satz} \label{sieb}
    Let $s \geq 2$ be a natural number and $\mathcal{P}$ a set of primes with the property $(s, p-1) \neq 1$ for all $p \in \mathcal{P}$. Denote the cardinality of this set by $P$. Furthermore, let $w(n)$ be a weight function subject to the above and additionally $w(n)=0$ for $n=0$ and $n \geq e^{P}$. Then
	\begin{align*}	
        S(\mathcal{A}) \ll &  P^{-1}\sum_{n=1}^{\infty} w(n) + P^{-2} \sum_{p \neq q \in \mathcal{P}} \left\lvert \sum_{n=1}^{\infty} w(n) \res sn{p} \overline{\res sn{q}} \right\lvert.
	\end{align*}
\end{satz}

\begin{bew}
	Let
 		$$ \Sigma = \sum_{n=1}^{\infty} w(n) \Bigg\lvert \sum_{p\in\sP} \res snp \Bigg\rvert^2.$$
	Whenever $n$ is an $s$-th power, $n=m^s$, say, for every $p$ coprime to $m$ the character takes the value $1$ and one has
		$$ \sum_{p \in \sP} \res snp = \dsum{p \in \sP} {p \nmid m} 1 \geq P - \sum_{p \mid m} 1.$$

	Now we split the last sum at $\log m$. The first term is
    	$$ \ll  \sum_{p < \log m} 1 \ll  \frac {\log m}{\log \log m}$$
 	by the Prime Number Theorem. Putting the number of prime divisors counted in the second sum equal to $r$, one observes that
 		$$ m > \prod_{\substack{ {p \geq \log m} \\ {p | m} }}p \geq (\log m)^r $$
    Hence the second term is $ =r < \frac{\log m}{\log \log m}$, as well. For $P < m < e^{P}$ this is $\ll \frac {P}{\log \log P}$, and for the remaining $m \leq P$  we have
		$$ \sum_{p \mid m} 1 \leq \sum_{p \leq m} 1 \leq \sum_{p \leq P} 1 \ll \frac{P}{\log P} $$
	by the Prime Number Theorem again. This establishes
		$$\sum_{p \in \mathcal{P} } \res snp \geq P + O \left( \frac{P}{\log \log P} \right) \gg P,$$
	so the entire expression is bounded by $\Sigma \gg P^2S(\mathcal{A})$.\\

	On the other hand, expanding the square yields
	\begin{align*}
		\Sigma =& \sum_{p, q \in \sP} \sum_{n=1}^{\infty} w(n) \res snp \overline{\res snq} \\
		\ll& \sum_{p \in \sP} \dsum{n=1}{ p \nmid n}^{\infty} w(n) + \sum_{p \neq q \in \sP} \sum_{n=1}^{\infty} w(n) \res snp \overline{\res snq} \\
		\ll& P \sum_{n=1}^{\infty} w(n) + \sum_{p \neq q \in \sP}  \left\lvert \sum_{n=1}^{\infty} w(n) \res snp \overline{\res snq} \right\lvert.
	\end{align*}
	The combination of these two results now proves the theorem.
\end{bew}

\begin{bems}
\begin{enumerate}[(a)]
    \item
        The condition $w(n) = 0$ for $n \geq e^P$ cannot be dispensed with. For illustration, let $m$ denote the product of the primes contained in a finite set $\mathcal P$, and define $w$ by $w(n_0) = 1$ for $n_0 = m^s$ and $w(n) = 0$ for any other $n \neq n_0$. In this case, the sieve gives $S(\mathcal{A}) \ll P^{-1}$, whereas $S(\mathcal{A}) = 1$ trivially.
    \item 	
        Theorem \ref{sieb} holds for arbitrary characters and also for appropriately normalized character sums, as long as they take the value $1$ on $s$-th powers. If, however, one chooses the principal character or a character sum it occurs in, independently from the choice of $\mathcal{P}$ and the weight function one finds
		  $$S(\mathcal{A}) \ll P^{-1}\sum_{n=1}^{\infty} w(n) + P^{-2} \sum_{p \neq q \in \mathcal{P}} \sum_{n=1} ^{\infty}w(n) \ll \sum_{n=1}^{\infty}w(n).$$
	   Since $S(\mathcal{A}) = \sum_{n=1}^{\infty}w(n^s) $, this is trivial.
    \item
        Theorem \ref{sieb} cannot be expected to be sharp even for non-principal characters. For instance, let $w(n) = 1 $ for $1 \leq n \leq x$, $w(n) = 0$  otherwise, and let $\mathcal{P}$ be the set of primes $p$ less than some $x^\alpha$ such that $s|p-1$ holds. Hence by the Siegel-Walfisz theorem one has $P \asymp x^\alpha (\log x)^{-1}$, and the number of $s$-th powers less than or equal to $x$ is given by
        	\begin{align*}
                S(\mathcal A)  & \ll  P^{-1} \sum_{n \leq x} 1 + P^{-2} \sum_{p \neq q \in \mathcal{P}} \left\lvert \sum_{n \leq x} \res snp \overline{\res snq}  \right\rvert \\
                & \ll  x^{1-\alpha} \log x + P^{-2} \sum_{p \neq q \in \mathcal{P}}   x^{\alpha} \log x  \\
                & \ll  x^{1-\alpha} \log x + x^{\alpha} \log x,
	       \end{align*}
        the second equation following from the P\'{o}lya-Vinogradov inequality. Thus, Theorem \ref{sieb} yields an estimate of at most $O(x^{1/2} \log x)$ for the number of $s$-th powers up to $x$, whereas obviously $S(\mathcal{A}) = x^{1/s} + O(1).$
\end{enumerate}
\end{bems}

\section{Application to twins of $s$-free numbers}

\subsection{Preliminaries}

In order to apply Theorem \ref{sieb} to the problem of twins of $s$-free numbers, we write
	$$E_s(n) =\sum_{j^s|n} \mu(j).$$
Then
	$$ \sum_{n \leq x} E_s(n) E_s(n+1) = \sum_{j,k} \mu(j) \mu(k) N(x,j,k), $$
where the counting function $N(x, j, k)$ is given by the number of solutions  of
	$$N(x, j, k) = \#\left\lbrace n \leq x: \; j^s | n, \; k^s | n+1\right\rbrace.$$
Note that $N(x, j, k) = x j^{-s} k^{-s} + O(1)$ for $(j,k) = 1$, whereas for $(j,k) > 1$ the condition is empty, since any common divisor of $j$ and $k$ would have to occur in $n$ as well as in $n+1$, which is impossible. The contribution of the terms with $jk\leq y$, where $y$ will be fixed later, is
\begin{align}
 	 x \dsum{jk \leq y}{(j,k)=1} \mu(j) \mu(k) (jk)^{-s} + O\Big( \sum_{jk \leq y} 1 \Big). \label{jk < y}
\end{align}

Now we complete the first sum to infinity. Assembling the $j$ and $k$ into $jk=n$, the different partitions of $n$ produce a factor $d(n)$. Therefore for \eqref{jk < y} we get
\begin{align*}
 	 = & x \sum_{(j,k)=1} \mu(jk) (jk)^{-s}  + O\Big( \sum_{n \leq y} d(n) \Big)+ O\Big(x \sum_{n>y} d(n) n^{-s} \Big) \\
 	 = & x \; \sum_{n} \frac {\mu(n) d(n)}{n^s} + O(y \log y)+ O(x y^{1-s+\epsilon})
\end{align*}
with the trivial bound $d(n) \ll n^{\epsilon}$ for the divisor function. The sum in the main term is over a multiplicative function and therefore can be rewritten as an Euler product:
	$$ \sum_{n} \frac {\mu(n) d(n)}{n^s} = \prod_p \sum_{i=0}^{\infty} \frac{ \mu(p^i)d(p^i) }{p^{is}} = \prod_p \left( 1 + \frac{-1 \cdot 2}{p^s} \right).$$
This is exactly the expression for the constant $C_s$ given in \eqref{cs}. Altogether, \eqref{jk < y} is equal to
  	$$ C_s x + O\left(y \log y\right) + O\left(xy^{1-s+\epsilon} \right).$$

The remaining $j$ and $k$ are contained in altogether $O\left( (\log x)^2 \right)$ intervals $J < j \leq 2J$ and $K < k \leq 2K$, where $JK \gg y$ and $J, K \ll x^{1/s}$. One can find $J$ and $K$ such that
	$$\sum_{jk>y} \mu(j) \mu(k) N(x, j, k) \ll N (\log x)^2,$$
where $N$ is given by
	$$N = \# \left\lbrace (j, k, u, v): \;  j^su+1=k^sv \leq x, \; J < j \leq 2J, \; K < k \leq 2K \right\rbrace. $$
Since we will choose $y$ in the range $x^{1/s} \leq y \leq x$, the resulting estimate is
\begin{align}
 	\sum_{n\leq x} E_s(n) E_s(n+1) = C_sx + O(y^{1+\epsilon}) + O\left( N (\log x)^2 \right). \label{2}
\end{align}

Without loss of generality, we can assume that $J \geq K$, where the sign of the $1$ in the expression for $N$ may change. For if $J<K$, swapping the parameters $j$ and $k$ as well as $u$ and $v$ transforms the condition $j^su +1 = k^sv$ into $j^su -1 = k^sv$ and one has $J \geq K$ as required.

Now in order to find a bound for $N$, we sort the quadruples counted by $N$ according to the value of $u$ and write $N = \sum N_u$, where
    $$ N_u = \# \left\lbrace (j, k, v): \; j^su \pm 1=k^sv \leq x, \; J < j \leq 2J, \; K < k \leq 2K \right\rbrace. $$

The $u$ are contained in $O(\log x)$ intervals $U < u \leq 2U$ with
\begin{align}
 	U \ll xJ^{-s}. \label{5}
\end{align}
A further examination of the equation $j^su \pm 1 = k^sv$ and the conditions on $j$ and $k$ allows us to write more conveniently
\begin{align}
	N_u  \leq & \# \left\lbrace (j, k, v): \; j^su \pm 1=k^sv, \; K < k \leq 2K, \; L \leq v \leq M\right\rbrace, \label{nu}
\end{align}
where the bounds on $v$ are given by
	$$L = \max\{2^{-s}K^{-s}(J^sU \pm 1), \;1\}, \quad M = K^{-s}(2^{s+1}J^sU \pm 1).$$

Putting
	$$N(U) = \sum_{U < u \leq 2U } N_u,$$
the estimate on $N$ is
\begin{align}
	 N \ll (\log x)  \max_{U \ll xJ^{-s} } N(U). \label{4}
\end{align}

\subsection{Application of Theorem \ref{sieb}}
In \eqref{nu} we found a representation of $N_u$ in which only the variables $k$ and $v$ on the right hand side of the equation are independently subject to size restrictions, so that all parameters on the left hand side are either fixed or uniquely determined by the equation.
This enables us to find a suitable weight function $w$ that translates the information contained in \eqref{nu} into the language of Theorem \ref{sieb} and thereby provides knowledge about the size of the $N_u$.

Let $w(n) = 0$ for any $n$ not being a multiple of $u^{s-1}$, and
	$$w(mu^{s-1}) = \# \left\lbrace (k, v): u|k^sv \mp 1, \; m=k^sv \mp 1, \; K<k \leq 2K, \; L \leq v \leq M\right\rbrace $$
otherwise.
This choice is motivated as follows: If $w(n) \neq 0$ is counted by $\mathcal{A}$, that is to say, if $n = mu^{s-1}$ is an $s$-th power, $m$ can be written as $m = j^su$. By construction, however, we also have $m = k^sv \mp 1$. Thus this representation translates pairs of $s$-free numbers into $s$-th powers, and the sieve of the previous section is applicable. Furthermore, by construction we have  $N_u \leq S(\mathcal{A})$.

Now let $\mathcal{P}$ denote the set of primes $p$ contained in the range $Q<p \leq 2Q$ subject to $p \nmid u$ and $(s, p-1) \neq 1$. Here, $Q$ lies in the interval
\begin{align}
 	(\log x)^2 \leq Q \leq x \label{6}
\end{align}
and will be chosen optimally later. By the theorem of Siegel-Walfisz we have $P \asymp Q(\log Q)^{-1}$ and therefore
    $$\log n \leq \log x \leq \sqrt{Q} \leq P$$
for $Q$ sufficiently large and for any $n$ counted with positive weight. Thus the numbers $n$ counted by $w$ are bounded above by $e^P$, as required, whereas the condition $w(0)=0$ is fulfilled trivially by construction. Theorem \ref{sieb} yields
\begin{align}\label{7}
    N_u \ll&  P^{-1} \sum_{n} w(n) + P^{-2} \negthickspace \sum_{p \neq q \in \sP} \left\lvert \sum_{k, v} \res s{u^{s-1}(k^sv \mp 1)}{p} \overline{\res s{u^{s-1}(k^sv \mp 1)}{q}} \right\rvert \nonumber \\
    \ll & \left( \frac {\log x}{Q} \right) \sum_{n} w(n) + P^{-2} \sum_{p \neq q \in \sP} \left\lvert \sum_{k, v} \; \res s{k^sv \mp 1}{p} \overline{\res s{k^sv \mp 1}{q}} \right\rvert,
\end{align}
where we used \eqref{6} and the multiplicity of characters. The restrictions on $k$ and $v$ in the inner sum are given by
\begin{align}
 	K<k \leq 2K, \quad L \leq v \leq M, \quad u|k^sv \mp 1, \label{8}
\end{align}
respectively. This expression has to be evaluated in order to find a bound on $N$ and therefore on
	$$\sum_{n \leq x} E_s(n)E_s(n+1) - C_sx.$$

The first term on the left hand side of \eqref{7} does not present any difficulties. The main effort will lie in finding an estimate for the second term. We will see that the sum splits up into independent factors that can be treated more easily. Here it will also become clear that the aspired saving is possible because the $s$-th powers are in some way well distributed modulo $pq$.

\subsection{Evaluation of the sifting}
The total contribution of the first term in \eqref{7} to $N(U)$ amounts to
\begin{align}
	\ll& \left( \frac{\log x}{Q}\right)  \sum_{k, v} \sum_{u|k^sv \mp 1} 1 \nonumber \\
	\ll& \left( \frac{\log x}{Q}\right)  \sum_{K < k \leq 2K} \; \dsum{n \leq x}{n \equiv \mp 1 \pmod{k^s}} d(n). \nonumber \\
	\ll& \left( \frac{\log x}{Q} \right) K^{1-s} \sum_{n\leq x}n^{\epsilon}  \nonumber\\
	\ll& Q^{-1} K^{1-s} x^{1+\epsilon}, \label{I}
\end{align}
where we used the standard estimate for the divisor function.\\

In order to understand the inner sum
	$$S := \sum_{k, v} \res s{k^sv \mp 1}{p} \overline{ \res s{k^sv \mp 1}{q} }$$
of the second term of \eqref{7} subject to the conditions \eqref{8}, one transforms the expression with the goal of separating the sums and encode any other information into exponential sums. Thereby we will obtain an estimate that will be easier to deal with analytically. Since additionally
	$$P^{-2} \sum_{p \neq q \in \mathcal{P}} 1 = O(1),$$
the bound on $\lvert S \rvert$ will be our final estimate on the second term in \eqref{7}.\\

Reducing the arguments of the characters modulo $upq$ and encoding the requirements on $k$ and $v$ into exponential sums, one obtains
\begin{align*}
    S = & \dsum{\alpha, \beta = 1}{u|\alpha^s \beta \mp 1}^{upq} \res s{\alpha^s\beta \mp 1}{p} \overline{ \res s{\alpha^s\beta \mp 1}{q} } \Big\lvert \dsum{K < k \leq 2K}{k \equiv \alpha \; (upq) } 1 \Big\rvert   \Big\lvert \dsum{L \leq v \leq M}{v \equiv \beta \; (upq) } 1 \Big\rvert\\
    = & \dsum{\alpha, \beta = 1}{u|\alpha^s \beta \mp1}^{upq} \res s{\alpha^s\beta \mp 1}{p} \overline{ \res s{\alpha^s\beta \mp 1}{q} } \left\lbrace \frac{1}{upq} \sum_{\gamma=1}^{upq} \ksum e\left( \frac{\gamma(\alpha-k)}{upq} \right) \right\rbrace  \\
    & \qquad \quad \times  \left\lbrace \frac{1}{upq} \sum_{\delta=1}^{upq} \vsum \e{\delta(\beta-v)}{upq} \right\rbrace.
\end{align*}
Now it is possible to separate the sums over $\alpha$ and $\beta$, over $k$, and over $v$, respectively. Writing
\begin{align}
    S(u, pq; \gamma, \delta)   = & \displaystyle{ \dsum{\alpha, \beta=1}{u|\alpha^s \beta \mp 1}^{upq} } \res s{\alpha^s\beta \mp 1}{p} \overline{ \res s{\alpha^s\beta \mp 1}{q} } \e{\gamma \alpha + \delta \beta}{upq}  \label{12} \\
    \vartheta_{\gamma}  = & \ksum \e{-\gamma k}{upq} \ll \min \left( K, \; \left\lVert \frac{\gamma}{upq} \right\rVert ^{-1} \right) \label{13} \\
    \varphi_{\delta}  = & \vsum \e{-\delta v}{upq} \ll \min \left( J^s K^{-s} U, \; \left\lVert \frac{\delta}{upq} \right\rVert ^{-1} \right), \label{14}
\end{align}
one has
\begin{align*}
 	S=&(upq)^{-2} \sum_{\gamma, \delta = 1} ^{upq} S(u, pq; \gamma, \delta)\vartheta_{\gamma} \varphi_{\delta}.
\end{align*}

By definition of the set $\mathcal{P}$, none of the $p \in \mathcal{P}$ divides $u$, so the sum $S(u, pq; \gamma, \delta)$ splits into factors, as shown in the lemma:
\begin{lem}
	Let $u = \prod r^f$ be the prime factor decomposition of $u$, and
	\begin{align*}
 		S_1(p; c, d) &= \sum_{\alpha, \beta = 1}^p \res s{\alpha^s\beta \mp 1}p \e{c\alpha + d\beta}{p} \\
 		S_2(r^f; c, d) &= \dsum{\alpha, \beta = 1}{r^f | \alpha^s\beta \mp 1}^{r^f} \e{c\alpha + d\beta}{r^f}.
	\end{align*}
	Let $p$ and $q$ be coprime to $u$. Then the factorization of the expression given in \eqref{12} reads
	\begin{align}
 		S(u, pq; \gamma, \delta) = S_1(p; c, d) \overline{S_1(q; -c, -d)} \prod_{\prod r^f = u} S_2(r^f; c, d), \label{label}
	\end{align}
	where $c$ and $d$ are integers such that  $(c, upq) = (\gamma, upq)$ and $(d, upq) = (\delta, upq)$.
\end{lem}
\begin{bew}
 	Firstly, for $p, q$ coprime, we observe that
		$$\sum_{a=1}^p \sum_{b=1}^q  \e{a}{p} \e{b}{q} = \sum_a \sum_b \e{aq+pb}{pq} = \sum_{c=1}^{pq} \e{c}{pq}$$
	by substituting $c=aq+bp$, where $c$ takes all values from $1$ to $pq$ exactly once. Thus, the product over the $S_2$ can be written as follows:
	   $$ \prod_{\prod r^f = u} S_2 (r^f; \gamma, \delta) = S_2(u; \gamma, \delta) = \dsum{\alpha, \beta = 1}{u| \alpha^s \beta \mp 1}^{u} \e{c\alpha + d\beta}{u}.$$
	This is not affected by the congruence condition on $\alpha$ and $\beta$.

	Furthermore, one has
	\begin{align*}
        \overline{S_1(q; -\gamma, -\delta)}  = & \sum_{\alpha, \beta = 1}^q \overline{\res {s}{\alpha^s \beta \mp 1}{p} } \; \; \overline{ \e{ - (\gamma \alpha + \delta \beta)}{q} } \\
 		= &  \sum_{\alpha, \beta = 1}^q \overline{\res {s}{\alpha^s \beta \mp 1}{p} } \e{\gamma \alpha + \delta \beta}{q}.
	\end{align*}

	Now, in order to calculate $S(u, pq; \gamma, \delta)$, put
	\begin{align*}
 		\alpha =& qu\alpha_1 + pu\alpha_2 + pq\alpha_3 \\
 		\beta = & qu\beta_1 + pu\beta_2 + pq\beta_3
	\end{align*}
	with
		$$ 1 \leq \alpha_1, \beta_1 \leq p \qquad 1 \leq \alpha_2, \beta_2 \leq q \qquad 1 \leq \alpha_3, \beta_3 \leq u. $$
	This is unique modulo the respective residues and one has
	\begin{align*}
        & \dsum{\alpha, \beta = 1 }{u| \alpha^s \beta \mp 1}^{upq} \res s{\alpha^s \beta \pm 1}p \overline{ \res s{\alpha^s \beta \pm 1}q } \e{\gamma \alpha + \delta \beta}{upq}\\
        =& \sum_{\alpha_1, \beta_1 = 1}^p \sum_{\alpha_2, \beta_2 = 1}^q \dsum{\alpha_3, \beta_3 = 1}{u | \alpha_3^s \beta_3 \mp 1}^{u} \res s{\alpha^s \beta \pm 1}p \overline{ \res s{\alpha^s \beta \pm 1}q } \\
        & \times \; \e{ (\gamma \alpha_1 + \delta \beta_1 ) qu + (\gamma \alpha_2 + \delta \beta_2) pu + (\gamma \alpha_3 + \delta \beta_3) pq }{pqu} \\
        =& \sum_{\alpha_1, \beta_1 = 1}^p \sum_{\alpha_2, \beta_2 = 1}^q \dsum{\alpha_3, \beta_3 = 1}{u | \alpha_3^s \beta_3 \mp 1}^{u} \res s{\alpha^s \beta \pm 1}p \overline{\res s{\alpha^s \beta \pm 1}q}\\
        & \times \; \e{\gamma \alpha_1 + \delta \beta_1}{p} \e{\gamma \alpha_2 + \delta \beta_2}{q} \e{\gamma \alpha_3 + \delta \beta_3}{u}.
	\end{align*}

	Since $p, q,$ and $u$ are relatively prime, the periodicity of the characters implies
	\begin{align*}
 		\res s{\alpha^s \beta \pm 1}p = & \res s{(qu\alpha_1 + pu\alpha_2 + pq\alpha_3 )^s (qu\beta_1 + pu\beta_2 + pq\beta_3) \pm 1}p \\
		=& \res s{(qu\alpha_1)^s (qu\beta_1) \pm 1}p.
	\end{align*}

	Altogether, we find
	\begin{align*}
        =& \sum_{\alpha_1, \beta_1 = 1}^p \res s{(qu\alpha_1)^s (qu\beta_1) \pm 1}p \e{\gamma \alpha_1 + \delta \beta_1}{p} \\
        &\times \sum_{\alpha_2, \beta_2 = 1}^q \overline{ \res s{(pu\alpha_2)^s (pu\beta_2) \pm 1}q } \e{\gamma \alpha_2 + \delta \beta_2}{q}  \dsum{\alpha_3, \beta_3 = 1}{u | \alpha_3^s \beta_3 \mp 1}^{u} \e{\gamma \alpha_3 + \delta \beta_3}{u} \\
        =& \sum_{\alpha_1, \beta_1 = 1}^p \res s{(qu\alpha_1)^s (qu\beta_1) \pm 1}p \e{qu(c \alpha_1 + d \beta_1)}{p} \\
        &\times   \sum_{\alpha_2, \beta_2 = 1}^q \overline{ \res s{(pu\alpha_2)^s (pu\beta_2) \pm 1}q } \e{pu(c \alpha_2 + d \beta_2)}{q} \\
        & \times \dsum{\alpha_3, \beta_3 = 1}{u | \alpha_3^s \beta_3 \mp 1}^{u} \e{pq(c \alpha_3 + d \beta_3)}{u},
	\end{align*}
	where $c$ and $d$ are determined uniquely modulo $upq$ by the congruences
	\begin{align*}
		\gamma &\equiv quc \pmod p; & \gamma &\equiv puc \pmod q; &  \gamma &\equiv pqc \pmod u\\
		\delta &\equiv qud \pmod p; & \delta &\equiv pud \pmod q; & \delta &\equiv pqd \pmod u.
	\end{align*}

	After shifting the summation indices $\alpha_1$ to $qu \alpha_1$ and correspondingly for $\alpha_2$,$\alpha_3$, the last expression reads
	\begin{align*}
        =& \sum_{\alpha_1, \beta_1 = 1}^p \res s{\alpha_1^s \beta_1 \pm 1}p \e{c \alpha_1 + d \beta_1}{p} \sum_{\alpha_2, \beta_2 = 1}^q \overline{ \res s{\alpha_2 ^s \beta_2 \pm 1}q } \e{c \alpha_2 + d \beta_2}{q} \\
        &\times  \dsum{\alpha_3, \beta_3 = 1}{u | \alpha_3^s \beta_3 \mp 1}^{u} \e{c \alpha_3 + d \beta_3}{u},
	\end{align*}
	which is $S_1(p; c, d) \overline{S_2(q; -c, -d)} S_2(u; c, d)$ as claimed.

    Since $p, q$ and $r$ are relatively prime, the gcd-conditions $(c, upq) = (\gamma, upq)$ and $(d, upq) = (\delta, upq)$ can be inferred directly from the determining equations for $c$ and $d$. Thus, the lemma is proved.
\end{bew}

\subsection{Exponential sums}

In order to deal with the exponential sums that have appeared here, we use the following lemma:
\begin{lem} \label{expo}
 	For the exponential sums $S_1$ and $S_2$ we have the following estimates:
	\begin{enumerate}[(i)]
		\item $\lvert S_1(p; c, d) \lvert \ll   p, \label{s1}$
		\item $\lvert S_2(p; c, d) \rvert \leq  s(s+1) \;p^{1/2}(p, c, d)^{1/2} + (s+1)^2 \label{s2}$.
	\end{enumerate}
\end{lem}

\begin{bew}
    First, we consider $S_1$. For $\alpha \neq p$  let $\bar{\alpha}$ be the multiplicative inverse of $\alpha$ modulo $p$. Shifting the summation index $\beta$ to $\bar{\alpha}^{s} \pm \beta$ yields
	\begin{align*}
 		\left \lvert S_1(p; c, d) \right \rvert =& \left\lvert \sum_{\alpha, \beta =1}^p \res s{\alpha^s\beta \mp 1}p \e{c\alpha + d\beta}{p} \right \rvert \\
 		\leq&  \left \lvert \sum_{\alpha = 1}^{p-1} \sum_{\beta=1}^p \res s{\beta}p \e{c\alpha + d \bar{\alpha} ^{s} \pm d\beta}{p} \right \rvert + O(p)\\
        =& \left \lvert \sum_{\alpha = 1}^{p-1} \e{c\alpha + d\bar{\alpha}^{s}}{p} \right \rvert \left \lvert \sum_{\beta=1}^p \res s{\beta}p \e{\pm d\beta}{p} \right \rvert + O(p) \\
 		=& \left \lvert S_2(p; c, d) \right \rvert  p^{1/2} + O(p) .
	\end{align*}
	Thus in the case $(p,c,d)=1$ the statement in \eqref{s1} follows directly from \eqref{s2} by using the well-known boundary on Gaussian sums.

    In the case that $p|(c,d)$ one has by the same strategy as before
    \begin{align*}
 		\left \lvert S_1(p; c, d) \right \rvert = \left\lvert \sum_{\alpha, \beta =1}^p \res s{\alpha^s\beta \mp 1}p \right \rvert
 		\leq&  \left \lvert \sum_{\alpha = 1}^{p-1} \sum_{\beta=1}^p \res s{\beta}p  \right \rvert + O(p) = O(p)
	\end{align*}
    and the statement follows directly.\\

	For estimating $S_2$ we distinguish four cases: The assertion is trivial if $p|(c, d)$. So it is if $p$ divides $d$, but not $c$, for
		$$S_2(p;c, d) = \sum_{\alpha = 1}^{p-1} \e{c\alpha}p = -1.$$
	In the opposite case that $p$ divides $c$, but not $d$, one has
		$$S_2(p; c, d) = \sum_{\beta = 1}^p \e{d\beta}p \dsum{\alpha=1}{\alpha^s\beta \equiv \pm 1 \pmod p}^p 1.$$
	For $\beta \neq p$, the congruence in the last summation has exactly $(s, p-1)$ solutions if $\pm \bar{\beta}$ is an $s$-th power; else it is insoluble.
	
    By choosing
        $$\dsum{\chi \bmod p}{\chi^s = \chi_0} \chi(n) = \begin{cases} \#\left\lbrace \chi \bmod p : \; \chi^s = \chi_0 \right\rbrace = (s, p-1) & n \hbox{ is } s \hbox{-th power} \\ 0 & \hbox{otherwise} \end{cases} $$
	as an indicator function on th $s$-th powers, changing the summation order and completing the sum over $\beta$ to a complete period, one finds
	\begin{align*}
        \left\lvert S_2(p; c, d) \right \rvert =& \Big \lvert \sum_{\beta = 1}^{p-1} \e{d\beta}p \dsum{\chi \bmod p}{\chi^s = \chi_0} \chi(\mp \bar{\beta}) \Big \rvert \\
 		=& \Big \lvert \dsum{\chi \bmod p}{\chi^s = \chi_0} \sum_{\beta = 1}^{p} \chi(\mp \bar{\beta}) \e{d\beta}p -1 \Big \rvert \\
		\leq& \dsum{\chi \bmod p}{\chi^s = \chi_0} p^{1/2}  + 1 \leq  (s, p-1) \; p^{1/2} + 1,
	\end{align*}
	where the estimate on Gaussian sums was used again. Since $(s, p-1) < s(s+1)$, the assertion follows.

    Finally, if neither $c$ nor $d$ is a multiple of $p$, we resort to a result from algebraic geometry that essentially goes back to Bombieri (see Chalk/Smith \cite{chalk}):
	\begin{satz}
	 	Let $\psi$ and $f$ be polynomials in $\mathbb{F}_p[X, Y]$ with  $\deg \psi = d_1$, $\deg f = d_2$.Under the condition that
		\begin{align*}
 			&f(X, Y) \not\equiv a \;(\hbox{mod} \; \psi_1(X, Y)) \; \hbox{in} \; \mathbb{F}_p \\
 			& \qquad \hbox{for all} \; a \in \mathbb{F}_p \hbox{ and for all absolutely irreducible } \psi_1|\psi \hbox{ in } \mathbb{F}_p,
		\end{align*}
		one has
			$$\Big\lvert \dsum{x, y \in \mathbb{F}_p}{\psi(x, y) = 0} \e{f(x, y)}{p} \Big\rvert \leq (d_1^2 - 3d_1 + 2d_1d_2)p^{1/2} + d_1^2.$$
	\end{satz}

    Letting $\psi(\alpha, \beta) = \alpha^s \beta \mp 1$ and $f(\alpha, \beta) = c\alpha + d\beta$, so $d_1 = s+1$ and $d_2 = 1$, the condition is fulfilled and the theorem yields
	\begin{align*}
 		\lvert S_2(p; c, d) \rvert \leq& \left( (s+1)^2 - 3(s+1) +2(s+1) \cdot 1 \right) \; p^{1/2} + (s+1)^2 \\
 		\leq &  s (s+1) \; p^{1/2} + (s+1)^2 \\
 		\leq& s(s+1) \; p^{1/2} + (s+1)^2,
	\end{align*}
	which is the desired result.
\end{bew}

\subsection{Evaluation continued}
Now that we can handle the exponential sums $S_1$ and $S_2$, we can continue in the evaluation of \eqref{label}. Let $w$ denote the product of all those prime factors $r|u$ that appear exactly once in the prime decomposition of $u$. Applying Lemma \ref{expo}\eqref{s2} produces a factor $s(s+1)$ for any prime occurring in $w$. In the product this amounts to
	$$(s(s+1))^{\nu(w)} = d_{s(s+1)}(w) \ll w^{\epsilon},$$
where $d_k(n)$ denotes the generalized divisor function as usual, and $\nu(n)$ is the number of prime divisors of $n$.

In the following we will not apply any diligence in distinguishing the $\epsilon$, but use the same symbol for any arbitrarily small exponents.

By Lemma \ref{expo} and the trivial estimate
	$$\lvert S_2(p^f; c, d) \rvert \leq p^f $$
for $f \geq 2$, the product of the $S_2$ is bounded by
\begin{align*}
 	\prod \left \lvert S_2(r^f; c, d) \right \rvert  \leq & U w^{-1} \left( w^{1/2 + \epsilon} (w,\gamma, \delta)^{1/2} + (s+1)^2 \right)  \\
 	\ll & U w^{-1/2 + \epsilon} (w, \gamma, \delta)^{1/2}.
\end{align*}

Altogether, this yields
	$$\lvert S(u, pq; \gamma, \delta) \rvert \ll Q^2 U w^{-1/2 + \epsilon} (w, \gamma, \delta)^{1/2}$$
and hence
\begin{align}
 	\lvert S \rvert \ll Q^{-2} U^{-1} w^{-1/2 + \epsilon} \sum_{\gamma, \delta}  \vartheta_{\gamma} \varphi_{\delta}  (w, \gamma, \delta)^{1/2}. \label{S}
\end{align}

Taking into account the bounds \eqref{13} and \eqref{14} as well as the symmetry of the distance function, one has further
\begin{align*}
 	&\sum_{\gamma, \delta}  \vartheta_{\gamma} \varphi_{\delta}  (w, \gamma, \delta)^{1/2} \\
 	\ll & J^s K^{1-s} U w^{1/2} + J^s K^{-s} Q^2 U^2 \sum_{1 \leq \gamma \leq \frac12upq} \gamma^{-1} (w, \gamma)^{1/2} \\
    & + \;K Q^2 U \sum_{1 \leq \delta \leq \frac12upq} \delta^{-1}(w,\delta)^{1/2} + Q^4 U^2 \sum_{1 \leq \gamma, \delta \leq \frac12upq} (\gamma \delta)^{-1} (w, \gamma, \delta)^{1/2},
\end{align*}
where the first term arises in the case that $\gamma \equiv \delta \equiv 0 \pmod{upq}$, the second and the third ones appear when exactly one of the two variables takes the value $0$, and the last one emerges in the case that neither $\gamma$ nor $\delta$ vanishes.

In order to evaluate the emerging sums, we consider as an example the first one over $\gamma$; the other ones behave analogously. One has
\begin{align*}
    \sum_{1 \leq \gamma \leq upq} (w,\gamma)^{1/2} \gamma^{-1}  \leq & \sum_{d|w} d^{1/2} \dsum{1 \leq \gamma \leq upq}{d|\gamma} \gamma^{-1} \ll  \sum_{d|w} d^{1/2} \dsum{e }{1 \leq de \leq upq} \frac{1}{de}\\
    \ll & \sum_{d|w} d^{-1/2} \log x \ll d(w) \log x \ll w^{\epsilon} \log x.
\end{align*}

Plugging these results into \eqref{S} yields
\begin{align*}
 	|S|  \ll & J^s K^{1-s} Q^{-2} w^{\epsilon} + \left\lbrace J^s K^{-s}U + K + Q^2 U \right\rbrace w^{-1/2+\epsilon} (\log x)^2.
\end{align*}

The contribution of $S$ to $N(U)$ can now be computed by summing over $u$. Firstly, one has
	$$\sum_{U < u \leq 2U} w^{\epsilon} \ll \sum_{U < u \leq 2U} u^{\epsilon} \ll U^{1+\epsilon}.$$
Secondly, $u$ decomposes uniquely into $u = wt$, where every prime factor occurs at least twice. The number of such square-full $t \leq z$ is easy to calculate: Any square-full $t$ can be written uniquely as $t = a^2b^3$ with (not necessarily coprime) integers $a$ and $b$. Thus, one has
	$$\sum_{a^2b^3 \leq z} 1 = \sum_{b \leq z^{1/3}} \left \lfloor \sqrt{ \frac{z}{b^3}} \right \rfloor  \leq z^{1/2} \sum_{b \leq z^{1/3}} b^{-3/2} $$
and the sum in the last term is $O(1)$.\\

For us, this implies
\begin{align*}
 	\sum_{U < u \leq 2U} w^{-1/2+\epsilon}  \leq &  \sum_{w \leq 2U} w^{-1/2 + \epsilon} \sum_{\substack{U < wt \leq 2U\\  t \; \mathrm{ squarefull}}} 1 \\
	\ll & \sum_{w \leq 2U} w^{-1/2 + \epsilon} (U/w)^{1/2} \\
	\ll & U^{1/2} \sum_{w \leq 2U} w^{-1+\epsilon} \ll U^{1/2 + \epsilon}.
\end{align*}

Altogether, taking into account \eqref{5}, the second term in \eqref{7} makes a contribution of
\begin{align}\label{II}
   	\ll  &   K^{1-s}J^sQ^{-2}U^{1+\epsilon} + J^{s}K^{-s}U^{\frac32+\epsilon} + KU^{\frac12+\epsilon} + U^{\frac32+\epsilon}Q^2\qquad \nonumber\\
  	\ll  &   x^{1+\epsilon}K^{1-s}Q^{-2} + x^{\frac32+\epsilon} J^{-\frac s2}K^{-s} + x^{\frac12 +\epsilon} J^{-\frac s2} K + x^{\frac32+\epsilon} J^{-\frac32s} Q^2.
\end{align}

\subsection{Conclusion of the proof}

Now we can combine the results from \eqref{I} and \eqref{II} and so, by choosing suitably the still undeterminated parameters $y$ and $Q$, find the requested estimate for
	$$\sum_{n \leq x} E_s(n) E_s(n+1) - C_sx.$$

First, one notes that $N(U)$ is bounded by
\begin{align*}
	\ll & K^{1-s} Q^{-1} x^{1+\epsilon} + x^{1+\epsilon}K^{1-s}Q^{-2} + x^{\frac32+\epsilon} J^{-\frac s2}K^{-s} + x^{\frac12 +\epsilon} J^{-\frac s2} K \\
	& + \; x^{\frac32+\epsilon} J^{-\frac32s} Q^2 \\
    \ll & x^{\frac32+\epsilon} J^{-\frac32 s} Q^2 + x^{\frac32+\epsilon} J^{-\frac s2} K^{-s} + x^{1+\epsilon} K^{1-s} Q^{-1}  + x^{\frac12+\epsilon} J^{-\frac s2} K.
\end{align*}

Comparing the first and the third term yields
	\begin{align*}
 		Q = x^{ - \frac{1}{6} + \epsilon} J^{ \frac s2 } K^{- \frac{s-1}{3}} + (\log x)^2;
	\end{align*}
because of $J \leq x^{1/s}$ this satisfies \eqref{6}. \\

With $J \geq K$ and $JK \gg y$, we now have
\begin{align*}
    N(U) \ll & x^{\frac32+\epsilon} J^{-\frac32s} (\log x)^2 + x^{\frac32+\epsilon} J^{-\frac s2} K^{-s} + x^{ \frac76 +\epsilon} J^{- \frac s2 } K^{ - \frac{2(s-1)}{3} } \\
    & + \; x^{\frac12+\epsilon}J^{-\frac s2}K \\
    \ll &  x^{\frac32+\epsilon} (JK)^{-\frac34s} +  x^{\frac32+\epsilon} J^{-\frac s2} K^{-s} + x^{ \frac76 +\epsilon} (JK)^{- \frac{7s-4}{12} } \\
    \ll &  x^{\frac32+\epsilon} y^{-\frac34s} +  x^{\frac32+\epsilon} J^{-\frac s2} K^{-s} + x^{ \frac76 +\epsilon} y^{- \frac{7s-4}{12} },
\end{align*}
whence, by \eqref{2} and \eqref{4}, we can conclude for the final estimate that
\begin{align*}
 	&\sum_{n \leq x} E_s(n) E_s(n+1) - C_sx  \\
 	\ll &  y^{1+\epsilon} + x^{\frac32+\epsilon} y^{-\frac34s} + x^{\frac32+\epsilon} J^{-\frac s2} K^{-s} + x^{ \frac 76 +\epsilon} y^{- \frac{7s-4}{12} }.
\end{align*}

If $y$ is chosen to equal $x^{ \frac{14}{7s+8} }$ , the first and the fourth term are of the same magnitude, therefore the expression above is
\begin{align}
 	\ll x^{ \frac{14}{7s+8} +\epsilon } + x^{\frac32 + \epsilon} J^{-\frac s2} K^{-s}. \label{J>K}
\end{align}
In the first term we recognize already the postulated error. \\

Now, in order to estimate the remaining term $x^{\frac23 + \epsilon}J^{- \frac s2}K^{-s}$, we need an elementary auxiliary bound for $N$:
\begin{lem} \label{hs}
 	For
		$$N = \# \left\lbrace (j, k, u, v): \; J < j \leq 2J, \; K < k \leq 2K, \;  j^su \pm 1=k^sv \leq x\right\rbrace $$
	we have the following estimate:
	\begin{align*}
 		N  \ll  x(J^{-s}K+(JK)^{1-s})(\log x)^{s-1}. \label{3}
	\end{align*}
\end{lem}
\begin{bew}
	We have
	\begin{eqnarray*}
 		N & = & \sum_{K < k \leq 2K} \; \sum_{u \leq xJ^{-s}} \dsum{J < j \leq 2J}{j^su \equiv \mp 1 \pmod{k^s}} 1.  \label{n1}
 	\end{eqnarray*}
    The congruence $j^su \equiv \mp 1 \pmod{k^s}$ can be evaluated by decomposing the modulus into prime powers and then applying Hensel's lemma on congruences modulo prime powers. Let $\xi$ be a solution of $j^su \pm 1 \equiv 0 \pmod{p^{t-1}}$ for a prime $p$ and $t \geq 2$. Then, according to Hensel's lemma, if
		$$\left. \frac{d}{dj} \right\rvert_{j = \xi} \left( j^su \pm 1 \right) = \left( s\xi^{s-1}u \right) \not\equiv 0 \pmod{p},$$
    any solution $\xi$ of $j^su \pm 1 \equiv 0 \pmod{p^{t-1}}$ can be mapped bijectively onto a solution  $\xi'$ of $j^su \pm 1 \equiv 0 \pmod{p^t}$. This is however the case, as $\xi$ is a solution to $j^su \pm 1 \equiv 0 \pmod p$  and therefore $u \equiv  \mp \xi^{-s} \pmod{p}$. Since by definition none of the $p \in \mathcal P$ divides $u$ (and therefore $\xi$) or $s$, one has
		$$ s\xi^{s-1}u \equiv  \mp s\xi^{s-1}\xi^{-s} \equiv \mp s\xi^{-1} \not\equiv 0 \pmod{p}.$$
    Thus the equivalences $j^su \pm 1 \equiv 0 \pmod{p^t}$ and $j^su \pm 1 \equiv 0 \pmod{p}$ have the same number of solutions, which, according to a theorem of Lagrange, is bounded above by the degree $s$ of the polynomial.
	
    Since the number of solutions modulo $k^s$ equals the product of the numbers of solutions of the single prime powers, the congruence $j^su \equiv \mp 1 \pmod{k^s}$ has at most $ s^{\nu(k)} \ll d_s(k)$ solutions. Altogether one has
	\begin{align*}
 		N &\ll \sum_{K<k\leq 2K} \; \sum_{u \leq xJ^{-s}} \min \left(d_s(k), \; JK^{-s} \right)   \\
		& \ll xJ^{-s}(1+JK^{-s})\sum_{K <k \leq 2K} d_s(k)  \\
		& \ll  x(J^{-s}K+(JK)^{1-s})(\log x)^{s-1}.
	\end{align*}
\end{bew}

Now, since $N$ is bounded by \eqref{J>K} as well as by the expression given in Lemma \ref{hs}, these two bounds can be combined in order to prove that all terms in $J$ and $K$ are dominated by the term $x^{ \frac{14}{7s+8} +\epsilon }$  we have already found. Obviously,
\begin{align*}
 	\min( xKJ^{-s}, x^{\frac32} J^{-\frac s2} K^{-s} ) & \leq  (xKJ^{-s})^{\lambda} (x^{\frac32}J^{-\frac s2}K^{-s})^{1-\lambda}
\end{align*}
for any $\lambda \in [0, 1]$. Choosing $\lambda = \frac{s}{2+3s}$ yields
\begin{align*}
    = & x^{\frac{s}{2+3s} + \frac32 \left( 1-\frac{s}{2+3s} \right) + \epsilon} (JK)^{ \frac{s}{2+3s} - s \left( 1-\frac{s}{2+3s} \right) } \\
 	\ll & x^{\frac{s}{2+3s} + \frac32 \left( 1-\frac{s}{2+3s} \right) + \epsilon } y^{ \frac{s}{2+3s} - s \left( 1-\frac{s}{2+3s} \right) } \\
    \ll & x^{\frac{s}{2+3s} + \frac32 \left( 1-\frac{s}{2+3s} \right) + \left( \frac{14}{7s+8} \right) \left( \frac{s}{2+3s} - s \left( 1-\frac{s}{2+3s} \right) \right)  + \epsilon} \\
	\ll & x^{ \frac{39s+24}{21s^2+38s+16} + \epsilon },
\end{align*}
where again $J \geq K$ and $JK \ll y$ were used. Since for any $s$ in question the exponent arising here is $ \frac{39s+24}{21s^2+38s+16} < \frac{14}{7s+8} $, the theorem follows.

\begin{bem}
    Even the auxiliary bound in Lemma \ref{hs} yields a better error term than Carlitz \cite{carl}: With little more effort one finds \eqref{2} with $O(y \log x)$ instead of $O(y^{1+\epsilon})$. Now, by applying $J \geq K$ and  $JK \gg y$, the lemma yields $N \ll xy^{\frac{1-s}{2} } (\log x)^{s-1}$, whence the choice $y = x^{\frac{2}{s+1}}$ gives Theorem \ref{satz} with an error term of $O \left( x^{\frac{2}{s+1}}(\log x)^3  \right)$. Our improvement comes mainly from being able to save a magnitude of $p^{1/2}$ over the trivial estimate in Lemma \ref{expo}\eqref{s2}, albeit at the cost of a constant $h(s)$. This evolves into the factor $d_{h(s)}(w) w^{-1/2}$ and finally, after the summation over $u$, produces a saving of $U^{1/2}$ at the cost of some $\log$ powers. The choices of $Q$ and $y$ then lead to the improvement of $O \big( x^{ \frac{1}{s^2} } \big)$ over the old result.
\end{bem}

\end{document}